\documentclass[12pt]{article}
\begin{document}
\author{Nguyen Van Chau\thanks{  Supported in part by the National
Basic 
program on Natural Science, Vietnam}\\{ \small Hanoi Institute of
Mathematics, 
P.O. Box 631, Boho 10000, Hanoi, Vietnam}\\ {\small E-mail: 
nvchau@thevinh.ncst.ac.vn }} 
\title{ Polynomial maps of complex plane with branched  valued set
isomorphic  
to complex line }

\date { May 10, 2001}

\maketitle

\begin{abstract}
We present a completed list of polynomial dominanting maps of $C^2 $ 
with 
branched  value curve isomorphic  to the complex line $ C $, up to
polynomial 
automorphisms.

{\it AMS Classification }: 14 E20, 14 B25.
\end{abstract}

\bigskip

{\bf 1.} 
Let $f:C^n \longrightarrow C^n$ be a polynomial dominanting map,  
$Close(f(C^n))=C^n $,  and denote by $\deg_f$ the geometric degree of $f$ - the
number of solutions of the  equation $f=a$ for generic points $ a\in
C^n$.  The 
{\it branched  value set } $E_f$ of  $f$  is defined  to be the
smallest subset of 
$C^n$ such that the map  
$$f  :C^n \setminus f^{-1}(E_f) \longrightarrow C^n\setminus E_f
\eqno (*)$$
gives a unbranched $\deg_f-$sheeted covering. It is well-known (see
[M]) that 
the branched value set $E_f$ is either empty set or an algebraic
hypersurface 
and  $ E_f=\{ a\in C^n : \# f^{-1}(a) \neq  \deg_f \}.$
If $ E_f=\emptyset $, then $f$ is injective, and hence, $f$ is an
automorphism of 
$C^n$ by the well-known fact that injective polynomial maps of $C^n$
are 
automorphisms (see [R]). The famous Jacobian conjecture ([BCW]) asserts
that 
$f$ must have a singularity if $E_f\neq \emptyset$.  The question
naturally 
raises as what can be said about polynomial dominanting maps of $C^n$
if their 
branched  value sets are isomorphic to a given algebraic hypersurface
$E$. 

In this article we consider polynomial dominanting maps of $C^2$ with 
branched  value set isomorphic to the complex line $ C $. We are
interesting in 
finding  a list of  such maps, up to polynomial automorphisms.  We say
that two 
map $f,g:C^2 \longrightarrow C^2$ are {\it equivalent } if there are
polynomial 
automorphisms $\alpha$ and $\beta $ of $C^2$ such that $\alpha\circ
f\circ 
\beta =g$.

\medskip

\noindent{\bf Theorem 1.} {\it A polynomial dominanting map $f$ of
$C^2$ 
with finite fibres and branched  value set isomorphic to $ C $ is
equivalent to  
the   map $ (x,y) \mapsto (x^{deg_f},y).$}

\medskip

In view of this theorem the equivalence classes of polynomial
dominanting maps 
of $C^2$ with finite fibres and branched  value set  isomorphic to
$ C $ are 
completed determined by the geometric degree of maps. Theorem 1 is an 
immediate consequence of  the following

\medskip

\noindent{\bf Theorem 2.} {\it A polynomial dominanting map of $C^2$ 
with branched  value set  isomorphic to $ C $ is equivalent to one of
the  
following maps

i) $ (x.y) \mapsto (x^{deg_f},y)$;

ii) $ (x,y)\mapsto (x^{deg_f},x^my)$, $ m \geq 1$.

iii) $ (x,y)\mapsto (x^{deg_f},x^m(x^ny +\sum_{i=0}^{n-1}a_ix^i))$,  $m
\geq 
1$, $n \geq 1$,  $a_0 \neq 0$ and $a_i=0$ for $i+m=0 (mod  \deg_f)$. }

\medskip

In fact, in the above list only maps of Type (i) have finite fibres.
The fiber at $ 
(0,0)$ of a map of the types (ii-iii) is the line $x=0$. Further,  as
shown  in \$ 4, 
topology of maps of Type (ii) and Type (iii) are quite different. The
proof of 
Theorem 2 presented in \$ 3  is an application of the famous theorem of

Abhyankar, Moh and Suziki ([AM],[S]) on the embedding of the complex
line 
into the complex plane. \$ 4 is devoted to some remarks and open
questions.

\medskip

{\bf 2.}
Let us recall some elementary facts on topology of polynomial in
two-variables. 
Let  $h(x, y) \in  C [x, y]$. By the {\it exceptional value set } $E_h$
of $h$  we 
mean a minimal set $E_h \subset  C $ such that the map 
$$h: C^2 \setminus h^{-1}(E_h) \longrightarrow  C  \setminus E_h$$
gives a smooth locally trivial fibration.  The set $E_h$ is at most a
finite set (see 
[V]).  The fiber  of this fibration,  denoted by $\Gamma_h$,  is called
{\it the 
generic fiber } of $h$. A polynomial $h$ is  {\it primitive } if its
generic fiber is 
connected. By the Stein factorization  a polynomial $h(x,y)$ can be
represented 
in the form $h(x,y)=\phi(r(x,y))$  for a primitive polynomial $r(x,y)$
and an one-
variable polynomial $\phi (t)$ (see [F]).

\medskip

The following lemma is an immediate consequence of Abhyankar-Moh-
Suzuki 
on the embedding of the complex line into the complex plane, which
asserts  
that  regular  embbedings of $ C $ in $C^2$ are equivalent to the
natural 
embbeding, or   equivalence, that  if $p(x,y)$ is irreducible
polynomial and if 
the curve $p=0$ is a smooth contractible  algebraic curve, then $p\circ
\alpha 
(x,y)=x$ for some polynomial automorphism $\alpha$ of  $C^2$.

\medskip

\noindent {\bf Lemma 1.} {\it Let $h \in  C [x,y]$. Suppose  that  the
generic 
fiber $\Gamma_h$ has $d$ connected components and each of them is 
diffeomorphic to $ C $. Then,  there exists a polynomial automorphism
$\alpha$ 
in $C^2$  such that }
$$h\circ \alpha (x,y)=x^d+a_1x^{d-1}+\dots +a_d.$$

\medskip
We will use this lemma in the situation when all of  fibres of $h$,
excepted at 
most one, are diffeomorphic to a dictinct  union of $d$ lines $ C $.
Then, the 
lemma shows that $h(\alpha (x,y))=x^d+a_d$ for an autommorphism
$\alpha$.

\medskip

\noindent {\it Proof of Lemma 1.} By the Stein factorization we can
represent $h
(x,y)=\phi(r(x,y))$ for a primitive polynomial $r\in  C  [x,y] $ and
$\phi \in  C  [t]
$. Further, one can choose $r$ and $\phi $ so that $\phi
(t)=t^{\deg\phi }+ $ {\it  
lower terms }. 

Observe that  for each $c\in  C $ the fiber $h^{-1}(c)$ consists of the
curves $r
(x,y)=c_i$, $i=1,\dots  \deg \phi $, where $c_i$ are zero points of
$\phi(t)-c=0$. 
Since the generic fiber $\Gamma_h$ has $d$ connected components and
each of 
them is diffeomorphic to $ C $,  the generic fiber $\Gamma_r$ of $r$ is 
diffeomorphic to $ C $, $\deg \phi =d$  and 
$$\phi (t) =t^d+a_1t^{d-1}+\dots +a_d.$$ 

Now, let $\gamma \in  C $ be a fixed generic value  of $r$. Then,  the
polynomial 
$r(x,y)-\gamma$ is irreducible and the curve $r(x,y)-\gamma=0$ is 
diffeomorphic to $ C $. So, in view of   Abhyankar-Moh-Suzuki Theorem 
there 
exists a polynomial automorphism $\alpha$ of $C^2$ such that $r(\alpha
(x,y))-
\gamma= x$. Then,  we get 
$$h(\alpha (x,y))=\phi (r(\alpha (x,y))=\phi (x+\gamma)=
x^d+a_1x^{d-1}+\dots 
+a_d.$$
Q.E.D\\  

\medskip

{\bf 3.} 
{\it Proof of Theorem 2:}
Let $f:C^2_{(x,y)} \longrightarrow C^2_{(u,v)} $ be a given
polynomial 
dominanting map with branched value set $E_f$ isomorphic  to $ C $,
where $ 
(x,y)$ and $ (u,v)$ stand for coordinates in $C^2$. 
In view of Abhyankar-Moh-Suszuki Theorem we can choose a polynomial 
automorphism $\alpha$ of $C^2$ so that the image $\alpha (E_f) $ is
the line 
$u=0$. 

Let $\bar f:=\alpha \circ f$, $\bar f=(\bar f_1,\bar f_2)$. Then,
$E_{\bar f}$ =
$\{u=0\}$, $\bar f^{-1}(E_{\bar f})=\bar f_1^{-1}(0)$
and the map  $$\bar f :C^2 \setminus \bar f_1^{-1}(0) \longrightarrow
C^2
\setminus \{u=0\}. $$ gives  a unbranched $\deg_f-$sheeted covering. 
This covering induces  unbranched $\deg_f-$sheeted coverings
$$\bar f :\bar f_1^{-1}(c) \longrightarrow \{u=c\}\simeq  C , \ c\neq
0.$$
Since $ C $ is simply connected,  
for every $0\neq c\in  C $  the fiber $\bar f_1^{-1}(c)$ consists of
exactly 
$\deg_f$ connected  components and  each of these components  is 
diffeomorphic to $ C $. So, by applying  Lemma 1 one can see that there
exists an 
automorphism $\beta $ of $C^2$ such that $$\bar f_1(\beta
(x,y))=x^{\deg_f}
+c.$$

Let  $\tilde f:=\bar f\circ \beta -(c,0)$. Then,  $\tilde
f(x,y)=(x^{\deg_f}, \tilde 
f_2(x,y))$.  Note that  $E_{\tilde f}=\{u=0\} $ and $\tilde
f^{-1}(\{u=0\})=\{x=0
\}$. So, by definition for each $ (a,b)\in C^2$, $a\neq 0$, the
equation $\tilde f 
(x,y)=(a,b)$ has exactly $\deg_f$ dictinct solutions. This implies that
for each $ 
(a,b)\in C^2$, $a\neq 0$,  the equation $\tilde f_2(\epsilon,y)=b$ has
a unique 
solution for each $\deg_f$ radical $\epsilon $ of $a$. Such a
polynomial $\tilde 
f_2(x,y)$  must be of the form 
$$\tilde f_2(x,y)= ax^ky + x^lg(x), 0\neq a \in  C , \ k \geq 0 \,\  l
\geq 0, \ g \in 
 C  [x].$$

For  $k=0$, we have
$\tilde f\circ \gamma (x,y)=(x^{\deg_f},y)$
for the automorphism $\gamma (x,y):=(x,a^{-1} y-ax^lg(x))$. 

Consider the case $k >0 $. Put  $m= \min \{ k;l\}$ and $n= k- m $.

For $ n=0 $ we can represent $\tilde f_2(x,y)$ in the form
$$\tilde f_2(x,y)=ax^m(y+ h(x))+c(x^{\deg_f}),$$
where $h(x),c(x)\in  C [x]$ and $$h(x)=\sum_{ i+m \neq 0 (mod
\deg_f)}a_ix^i.
$$
Put  

$\gamma_1 (u,v):=(u,a^{-1}(v-c(u)));$ 

$\gamma_2 (x,y):=(x,y-h(x)).$

\noindent Then, we get
$$\gamma_1\circ \tilde f\circ \gamma_2 (x,y) = (x^{\deg_f},x^my).$$

For the case $n >0$ we can represent
$$
\tilde f_2 (x,y)=ax^m(x^ny+ h(x)+x^n b(x))+c(x^{\deg_f}),
$$
where  $h(x), b(x),\ c(x)\in  C [x]$ and 
$$h(x)=\sum_{i=0, \dots n-1, \ i+m \neq 0 (mod \deg_f)}a_ix^i.$$
Put 

$\gamma_1 (u,v):=(u,a^{-1}(v-c(u)));$

$\gamma_2 (x,y):= (x,y-b(x)).$

\noindent Then, 
$$\gamma_1\circ \tilde f\circ \gamma_2
(x,y)=(x^{\deg_f},x^m(x^ny+h(x))).$$

Thus, the considered map $f$ is allway equivalent to one of maps of the
types (i-
iii).
Q.E.D\\ 

\medskip

{\bf 4.} Let us to conclude the paper by some remarks and open
questions.

\medskip

i) Topologically, the behavior of maps of Type (i), Type (ii) and Type
(iii) are 
quite different. The maps of Type (i) have finite fibres, while the
fiber at $ (0,0)$ 
of a map of the types (ii) or (iii) is the line $\{ x=0\}$.
Furthermore, for an 
irreducible germ curve $\gamma \subset C^2$  located at $ (0,0)$  and 
intersecting with  the line $\{u=0\}$ at $\{(0,0)\}$,  the inverse
image of 
$\gamma$  by a map of Type (ii) is connected and consists of the line
$x=0$ and 
a branch located at $ (0,0)$. But, the inverse image of $\gamma$ by a
map of 
Type (iii) consists of the line $x=0$ and a branch at infinity.

\medskip

ii) As shown in Theorem 2 and its proof, the polynomial dominanting
maps of 
$C^2$ with branched value curve isomorphic to $ C $ form a narrow
classes 
among polynomial dominanting maps of $C^2$.  The structure of  
covering (*) 
associated to  such  a map is very simple. Up to automorphisms of
$C^2$,  this 
covering looks like a unbranched covering from $C^2 \setminus \{ x=0\}
\longrightarrow C^2 \setminus \{ u=0\}$. The singular set of such a
map is 
isomorphic to $ C $. In particular,{\it  a polynomial map of $C^2$ with
branched 
value curve isomorphic to $ C $ must has a singularity }. This is true
even when  
the branched value curve is only homeomorphic to $ C $ (see in [C1]). It
is worth 
to present here the following  observation, which is reduced from the
results in 
[C2]. By a {\it $J$-curve } we mean a curve $E$ in $C^2$ of the form
$E=\cup_i \phi_i( C )$ for some polynomial maps $\phi_i: C 
\longrightarrow 
C^2$, $i=1, \dots k$, such that 

a) Each component  $\phi_i( C )$ is not isomorphic to $ C $, and

b) For every polynomial automorphism $\alpha$ of $C^2$
$$ {\deg p_1 \over \deg q_1} =... = {\deg p_k  \over \deg q_k}, $$
where $ \alpha \circ \phi_i:=(p_i,q_i). $

\medskip

\noindent {\bf Theorem 3.} (See [Thm. 4.4 , C2].) {\it  A polynomial 
dominanting map $f$ of $C^2$ with $E_f\neq \emptyset$ must has a 
singularity if $E_f$ is not  a $J$-curve .}

\medskip

\noindent In fact, it is shown in [C2]  that if $f$ has non-zero
constant jacobian 
and  if $E_f\neq \emptyset$, then $E_f$ is a $J$-curve  in above and
for every polynomial automorphism $\alpha$ of $C^2$
$$ {\deg p_i \over \deg q_i}={\deg (\alpha \circ f)_1  \over \deg (
\alpha \circ f)
_2}, \ i=1,\dots k,  $$
where $ \alpha \circ \phi_i:=(p_i,q_i) $ and $\alpha \circ f:=((\alpha
\circ f)_1,( 
\alpha \circ f)_2)$. In view of this theorem, the Jacobian  conjecture
in $C^2$ 
can be reduced to the question {\it whether  a $J$-curve could not be 
the 
branched value set of a nonsingular polynomial map of $C^2$. }  

Note that the situation is quite different in the case of holomorphic
maps. 
Orevkov in [O]  had constructed  a nonsingular holomorphic map from a 
Stein 
manifold homeomorphic to $ R^4$ onto an open ball in the complex plane
$C^2
$  which gives a  three-sheeted branching covering with the branched
value set 
diffeomorphic to $ R^2$. 

\medskip

iii) It is worth to  find analogous statements of Theorem 1 for
high-dimensional 
cases. Let $F=(F_1,F_2,\dots F_n): C^n_{(x_1,\dots x_n)
}\longrightarrow C^n_
{(u_1,\dots u_n) }$ be a polynomial dominanting map with finite fibres.

\medskip

{\bf Problem:} {\it Suppose $E_F=\{u_1=0\}$.  

i) Does there exists an automorphism $\alpha$  of $C^n$ such that 
$$F_1\circ \alpha (x_1,\dots x_n)=x_1^{\deg_F}\ ? $$

ii) Is $F$ equivalent to the map 
$$ (x_1,\dots x_n) \mapsto (x_1^{\deg_F},x_2,\dots x_n) \  ?$$
}

Such a map $F$ gives a locally trivial  fibration $F_1:C^n \setminus \
F^{-1}(0)
\longrightarrow  C  \setminus \{0\}$ with fiber diffeomorphic to a
distinct union 
of $\deg_f$ space $C^{n-1}$. Further, every connected component of the
fibers 
$F_1^{-1}(c), \ c\neq 0$, is isomorphic to $C^{n-1}$. The situation
here seems 
to be more simple than those in the problem of embeddings $C^{n-1}$
into 
$C^n$ - a generalization of the Abhyankar-Moh-Suzuki Theorem, which 
asks 
whether a regular  embbeding of $C^{n-1}$ in $C^n$ is equivalent to
the 
natural embbeding, which is still open for $n >2$.

\medskip

\noindent {\it Acknownlegments}: 
The author wishes to thank very much 
Prof. V.H. Ha for his valuable comments and suggestions.

\medskip 

\noindent {\bf References} 

\noindent [AM] S. S. Abhyankar and T. T.  Moh,    {\it  Embeddings of
the line 
in the 
plane},    J.  Reine Angew.  Math.  276 (1975),    148-166.  

\noindent [BCW] H.  Bass,    E.  Connell and D.  Wright,   {\it  The
Jacobian 
conjecture: 
reduction of degree and formal expansion of the inverse},    Bull. 
Amer.  Math.  
Soc.  (N. S.) 7 (1982),    287-330.  

\noindent [C1] Nguyen Van Chau,     {\it Remark on the Vituskin's
covering},   
Act.  Math.   Vietnam ,  24, No 1, 1999 , 109-115.  

\noindent [C2] Nguyen Van Chau, {\it 
Non-zero constant Jacobian polynomial maps of $ C^2$.}
J. Ann. Pol. Math. 71, No.3, 287-310 (1999). 

\noindent [F]  M. Furushima, {\it Finite groups of polynomial
automorphisms in 
the complex affine plane}, 1, Mem. Fac. Sci., Kyushu Univ., Ser. A, 36
(1982), 
83-105.

\noindent [M]  D. Mumford, {\it Algebraic geometry, I. Complex
projective 
varieties, } Spring-Verlagm Berlin Heidelberg New York, 1976.

\noindent [O] S.  Yu.   Orevkov,    {\it Rudolph diagram and analytical
realization of 
Vitushkin's covering},    Math.   Zametki,    60 (1996),    no 2,   
206-224,    319. 

\noindent [R] W. Rudin, {\it Injective polynomial maps are
automorphisms }, 
Amer. Math. Monthly 102 (1995), No 6, 540-543.

\noindent [S] M.  Suzuki,    {\it Propietes topologiques des polynomes
de deux 
variables compleces at automorphismes algebriques  de lespace $C^2$,  
} J.  
Math.  Soc.  Japan 26,    2 (1976),    241-257. 

\noindent [V] J. L Verdier,    {\it Statifications de Whitney et
theoreme de 
Bertini-Sard,   } Invent.  Math.  36 (1976),    295-312. 

\end{document}